\documentclass{amsart}%
\usepackage{amsmath}
\usepackage{graphicx}%
\usepackage{amsfonts}%
\usepackage{amssymb}
\newtheorem{theorem}{Theorem}[section]
\theoremstyle{plain}
\newtheorem{corollary}[theorem]{Corollary}
\newtheorem{lemma}[theorem]{Lemma}
\newtheorem{proposition}[theorem]{Proposition}
\theoremstyle{remark}
\newtheorem{remark}[theorem]{Remark}
\newtheorem{example}[theorem]{Example}
\newtheorem{definition}[theorem]{Definition}
\newtheorem{notation}[theorem]{Notation}
\numberwithin{equation}{section}

\begin{document}
\title[Isophasal Scattering Manifolds]{Continuous Families of Isophasal Scattering Manifolds}
\author{Carolyn Gordon}
\address[Gordon]{Department of Mathematics, Dartmouth College, Hanover, New Hampshire
03755, U. S. A.}
\email{Carolyn.S.Gordon@dartmouth.edu}
\author{Peter Perry}
\address[Perry]{Department of Mathematics, University of Kentucky, Lexington, Kentucky
40506--0027, U. S. A. }
\email{perry@ms.uky.edu.}
\thanks{Gordon supported in part by NSF\ grant DMS-0072534.}
\thanks{Perry supported in part by NSF\ grants DMS-9707051 and DMS-0100829.}
\date{July 24,2002}

\begin{abstract}
We construct continuous families of scattering manifolds with the same
scattering phase. The manifolds are compactly supported metric perturbations
of Euclidean $\mathbf{R}^{n}$ for $n\geq8$. The metric perturbation may have
arbitrarily small support.

\end{abstract}
\maketitle

\section{Introduction}

\label{sec.intro}

Inverse spectral geometry for compact Riemannian manifolds is the study of
what geometric properties of the manifold are determined by the eigenvalues of
the Laplacian. For non-compact Riemannian manifolds, there may be only
finitely many $L^{2}$-eigenvalues of the Laplacian (or even no $L^{2}%
$-eigevalues at all!)\ but there are several possible analogues of spectral
data for which one can pose a similar inverse problem. A particularly
attractive setting in which to study the inverse problem is Euclidean
$\mathbf{R}^{n}$ with a compactly supported perturbation of the Euclidean
metric: in what follows we will write $X=(\mathbf{R}^{n},g)$ where $g$ is such
a compactly supported metric perturbation. Our goal is to construct continuous
families of such manifolds with the same `spectral' data, appropriately defined.

For the class of manifolds that we will consider, the Laplacian has purely
continuous spectrum in $[0,\infty)$ and no $L^{2}$-eigenvalues. Thus, the
resolvent of the Laplacian is an analytic function $R(z)=(\Delta_{X}-z)^{-1}%
$on $\mathbf{C}\backslash\lbrack0,\infty)$; as we discuss in what follows, the
resolvent admits a meromorphic continuation to a double covering of the
complex plane if $n$ is odd, and a logarithmic covering of the complex plane
if $n$ is even (this result follows, for example, from the ``black box
scattering'' formalism introduced by Sj\"{o}strand and Zworski in
\cite{SZ:1991} [but this only treats $n$ odd]).  \emph{Resolvent resonances}
are poles of the meromorphically continued resolvent; they serve as discrete
data analogous to the eigenvalues but are less easily studied than the
eigenvalues since their presence signals the solution to a non-selfadjoint
eigenvalue problem for the underlying differential operator. In the literature
they are also referred to simply as resonances. For the case considered here,
the resolvent resonances are identical to the poles of the meromorphically
continued scattering operator, which are called \emph{scattering
resonances} (we discuss the scattering operator in greater detail in what
follows). We will call two such manifolds \emph{isopolar} if they have
the same scattering resonances with multiplicities.

For certain classes of non-compact manifolds, including the class to be
studied here, one can define the \emph{scattering phase}, a function roughly
analogous to the counting function for eigenvalues in the compact problem. We
define and discuss the scattering phase in section \ref{sec.scattering} of
what follows. Two noncompact Riemannian manifolds are said to be
\emph{isophasal} if they have the same scattering phase. In our case, if two
manifolds have the same scattering phase, they also have the same scattering
poles, so isophasality is a stronger condition than isopolarity.

While the past two decades have seen an explosion of examples of isospectral
compact Riemannian manifolds, there are relatively few examples known of
isopolar or isophasal manifolds. In dimension greater than one, the known
examples include finite-area Riemann surfaces (both isopolar and
isophasal--see Berard \cite{Berard:1992} and Zelditch \cite{Zelditch:1992a},
\cite{Zelditch:1992}), Riemann surfaces of infinite area (isopolar and
isophasal--see Guillop\'{e}-Zworski \cite{GZ:1997} and Brooks-Davidovich
\cite{BD:2002}), three-dimensional Schottky manifolds (isopolar--see
Brooks-Gornet-Perry \cite{BGP:1999}), and surfaces that are isometric to
Euclidean space outside a compact set (isopolar and isophasal--see
Brooks-Perry \cite{BP:2001}). In all these examples, the manifolds share a
common Riemannian covering. They are constructed by the analog of a technique
of T. Sunada \cite{Sunada:1985}, which produces compact isospectral manifolds
with a common finite covering.

We will prove:

\begin{theorem}
\label{thm.main}For every $n\geq8$, there exist continuous families of
isophasal, non-isometric Riemannian metrics on $\mathbf{R}^{n}$ which are
Euclidean outside of a compact set of arbitrarily small volume. There also
exist pairs of such metrics on $\mathbf{R}^{6}$.
\end{theorem}

\begin{remark}
Letting $m=n-4$, the parameter space for the continuous families of isophasal
metrics on $\mathbf{R}^{n}$ that we will construct has dimension
\[
d\geq\frac{m(m-1)}{2}-\left[  \frac{m}{2}\right]  \left(  \left[  \frac{m}
{2}\right]  +2\right)  >1
\]
if $n=9$ or $n\geq11$. (If $n=8$ or $n=10$, the parameter space has dimension
at least $1$).
\end{remark}

In addition to proving the existence of the continuous families, we will give
an explicit example of a triple of isophasal metrics on $\mathbf{R}^{12}$. We
will see that these metrics have very different geometry. Indeed their
isometry groups have different dimension and structure.

To our knowledge, the isophasal metrics of Theorem~\ref{thm.main} differ from the
other known examples of isophasal or isopolar metrics in the following ways:

\begin{itemize}
\item {} They are the first \emph{continuous} families of isophasal or isopolar metrics;

\item {} They are the first examples for which
the manifolds do not share a common Riemannian cover;

\item {} They are the first isophasal or isopolar compact metric perturbations of the
Euclidean metric on $\mathbf{R}^{n}$.
\end{itemize}

Just as the examples cited above were based on an extension to noncompact
manifolds of a technique first developed by Sunada for constructing
isospectral compact manifolds, our examples use an extension of a technique
involving torus actions previously developed for the construction of
isospectral compact manifolds with different local geometry. In fact the
metrics that we use here were first constructed in Gordon \cite{G:2001} and
Schueth \cite{Schueth:2002}, where they were restricted to balls and spheres.
The method of torus actions was used to show that these metrics on balls and
spheres were isospectral.

One might worry that the examples constructed have trivial scattering (e.g.,
have no scattering poles!). We show, however, that the isosphasal metrics can
always be chosen to have infinitely many resonances. For metric scattering on
$\mathbf{R}^{n}$, S\'{a} Barreto and Tang \cite{SbT:1998} ($n$ odd) and Tang
\cite{Tang:2000} ($n$ even) proved the existence of infinitely many resonances
so long as the second relative heat invariant $a_{2}$ is non-vanishing. They
also gave various geometric hypotheses which guarantee the non-vanishing of
$a_{2}$: one of these is that the given metric is not flat but is a compactly
supported perturbation of the Euclidean metric that is close in $\mathcal{C}%
^{k}$ topology to the Euclidean metric for sufficiently large $k$ (Theorem 1.3
of \cite{SbT:1998} and Theorem 1.1 of \cite{Tang:2000}%
)\footnote{\label{Footnote:Kuwabara}Both of these papers rely on a result of
Kuwabara \cite{K:1980} which states that, given a flat metric $\gamma$ on a
compact manifold $X$, there is a neighborhood $U$ of $\gamma$ in the
$\mathcal{C}^{\infty}$ topology so that if $a_{2}(g)=0$ and $g\in U$, then $g$
is flat. In fact, a close examination of \cite{K:1980} shows that it is
sufficient for $g$ and $\gamma$ to be close in $\mathcal{C}^{k}$ topology for
$k$ sufficiently large: see Theorem A' of section 6 in \cite{K:1980}. For the
connection between Kuwabara's result on compact manifolds and the result on
metric perturbations of Euclidean $\mathbf{R}^{n}$, see the proof of Theorem
1.3 of \cite{SbT:1998} which uses finite propagation speed for solutions of
the wave equation.}. In our examples, it is easily verified that the metrics
are not flat, and it is easy to construct examples where the metrics are
arbitrarily close in $\mathcal{C}^{k}$ sense to the Euclidean metric for any
large fixed $k$\footnote{As explained in Section~\ref{sec.examples}, the metrics depend on
a skew-symmetric bilinear form and a $\mathcal{C}_{0}^{\infty}(\mathbf{R}%
^{n})$ function $\varphi$ which defines the support of the perturbation. One
takes the function $\varphi$ sufficiently small in $\mathcal{C}^{k}$-sense.}.
We can actually remove the assumption that our metrics are $\mathcal{C}^{k}$
close to the Euclidean metric by computing the $a_{2}$ heat invariant
directly, at the cost of imposing a genericity assumption \footnote{The genericity
condition we impose is merely a genericity condition on the choice of cut-off function $\varphi$.} on the space of
metrics; we carry out this computation in Section \ref{sec.a2} for the metrics
on $\mathbf{R}^{n}$ with $n\geq9$. 

The plan of this paper is as follows. In section \ref{sec.scattering}, we
discuss basics of scattering theory for asymptotically Euclidean manifolds. In
section \ref{sec.technique}, we develop a method, based on the use of torus
actions, for constructing manifolds with the same scattering phase (see
Theorem \ref{thm.gen.method}). In section \ref{sec.examples} we apply the
technique of the previous section to show that the metrics on $\mathbf{R}^{n}$
constructed in \cite{G:2001} and \cite{Schueth:1999} are isophasal, thus
proving Theorem \ref{thm.main}. We give full details of the examples in
dimension $n\geq9$, based on the examples in \cite{G:2001} (modified as in
\cite{Schueth:1999} so that the metrics are Euclidean outside of an
arbitrarily small compact set). The lower-dimensional examples are given by
metrics constructed in \cite{Schueth:2002}. Since the methods of section 3
apply in exactly the same way to these examples, we do not include the details
here. Finally, in Section \ref{sec.a2}, we carry out the explicit computation
that the $a_{2}$ heat invariant is generically nonvanishing for our examples
in dimension 9 and above.

In a second paper, in preparation, we will construct pairs of
\emph{conformally equivalent} isosphasal Riemannian metrics, again equal to
the Euclidean metric outside of a compact set. We will also construct pairs of
isosphasal potentials for the Schr\"odinger operator on $(\mathbf{R}^{n},g)$,
where again $g$ is a compact perturbation of the Euclidean metric.

\emph{Acknowledgements}. Peter Perry thanks Carolyn Gordon and David Webb for
hospitality at Dartmouth College during part of the time that this work was
done. The authors acknowledge the support of the NSF\ under grant DMS-0207125
which supported a conference on Inverse Spectral Geometry at the University of
Kentucky during which this work was completed.

\section{Metric Scattering on $\mathbf{R}^{n}$}

\label{sec.scattering}

In this section we review scattering theory for manifolds $X=(\mathbf{R}%
^{n},g)$ where $g$ is a compactly supported metric perturbation of the
Euclidean metric: see especially \cite{Melrose:1992} and see
\cite{Melrose:1995} for an expository treatment that includes the case
considered here. Letting $\Delta_{X}$ be the positive Laplace-Beltrami
operator on $X$, it follows from the classical Rellich uniqueness theorem that
$\Delta_{X}$ has no $L^{2}$-eigenvalues, and it is easy to prove that
$\Delta_{X}$ has purely absolutely continuous spectrum in $[0,\infty)$. Thus
the resolvent operator $\widetilde{R}(z)=(\Delta_{X}-z)^{-1}$, considered as a
mapping from $L^{2}(X)$ to itself, is an operator-valued analytic function of
$z$ in $\mathbf{C}\backslash\lbrack0,\infty)$. It can be shown that the
mapping $R(\lambda)=\widetilde{R}(\lambda^{2})$, initially defined on the
half-plane $\Im(\lambda)>0$ and viewed as a map from $\mathcal{C}_{0}^{\infty
}(\mathbf{R}^{n})$ to $\mathcal{C}^{\infty}(\mathbf{R}^{n})$, admits a
meromorphic continuation to the complex $\lambda$-plane if $n$ is odd, and to
the logarithmic plane if $n$ is even. At any poles $\zeta$, the resolvent
admits a Laurent expansion with finite polar part of the form%
\[
\sum_{j=1}^{N_{\zeta}}\frac{A_{j}}{\lambda-\zeta}%
\]
where the $A_{j}$ are finite-rank operators from $\mathcal{C}_{0}^{\infty
}(\mathbf{R}^{n})$ to $\mathcal{C}^{\infty}(\mathbf{R}^{n})$. The multiplicity
of the pole $\zeta$ is defined as $\dim\left(  \oplus_{j}(\operatorname*{Ran}%
A_{j})\right)  $. 

To define the scattering phase, we first recall that the absolutely continuous
spectrum is parameterized by scattering solutions to the eigenvalue equation
$(\Delta_{X}-\lambda^{2})u=0$ which are easily described. In what follows,
write $x\in\mathbf{R}^{n}$ as $x=r\omega$ where $r\geq0$ and $\omega\in
S^{n-1}$.

\begin{proposition}
\label{prop.unique}Fix $f_{-}\in\mathcal{C}^{\infty}(S^{n-1})$ and $\lambda
>0$. There exists a unique solution of the equation
\[
(\Delta_{X}-\lambda^{2})u=0
\]
having the asymptotic form
\begin{equation}
u(r\omega)=r^{(1-n)/2}e^{i\lambda r}\,f_{+}(\omega)+r^{(1-n)/2}e^{-i\lambda
r}f_{-}(\omega)+\mathcal{O}(r^{-(n+1)/2}) \label{eq.scatt.asy}%
\end{equation}
as $r\rightarrow\infty$. In particular, the function $f_{+}\in\mathcal{C}
^{\infty}(S^{n-1})$ is uniquely determined.
\end{proposition}

For a proof see \cite{Melrose:1992}.

The Proposition implies that the mapping $f_{-}\mapsto f_{+}$ is a
well-defined mapping from $\mathcal{C}^{\infty}(S^{n-1})$ to itself. We denote
this map, the absolute scattering matrix for $X$, by $S(\lambda)$. >From the
definition, it is clear that $S(\lambda)$ is a linear mapping, and that
$S(\lambda)^{-1}=S(-\lambda)$ for real $\lambda\neq0$. In the case of
$X_{0}=(\mathbf{R}^{n},g_{0})$ (where $g_{0}$ is the Euclidean metric on
$\mathbf{R}^{n}$), we have the explicit formula
\[
u(x)=\int_{S^{2}}\exp(-i\lambda x\cdot\omega)\,f_{-}(\omega)\,d\omega
\]
and a stationary phase calculation shows that the absolute scattering matrix
$S_{0}(\lambda)$ is given by
\begin{equation}
\left(  S_{0}(\lambda)\varphi\right)  (\omega) =i^{n-1}\varphi(-\omega
)\text{.} \label{eq.s.free}%
\end{equation}
Since $X$ is a compactly supported metric perturbation of $X_{0}$, it is not
surprising that the `relative scattering matrix'
\begin{equation}
S_{r}(\lambda)=S(\lambda)S_{0}(\lambda)^{-1} \label{eq.s.rel}%
\end{equation}
has especially nice properties (see, for example, section 5.2 of
\cite{Melrose:1995}):

\begin{proposition}
For real $\lambda\neq0$, the relative scattering matrix $S_{r}(\lambda)$
extends to a unitary operator from $L^{2}(S^{n-1})$ to itself. Moreover
\[
S_{r}(\lambda)=I+T(\lambda)
\]
where $T(\lambda)$ is an integral operator with integral kernel belonging to
$\mathcal{C}^{\infty}(S^{n-1}\times S^{n-1})$.
\end{proposition}

In particular, $T(\lambda)$ extends to a trace-class operator on
$L^{2}(S^{n-1})$, so that the operator determinant
\[
\det S_{r}(\lambda)=\det(I+T(\lambda))
\]
is well-defined (see, for example, \cite{Simon:1979} for discussion of
operator determinants). Since$\,S_{r}(\lambda)$ is unitary, it follows that
$\det S_{r}(\lambda)$ has modulus one. We note for use later that if $A$ is a
trace-class operator on a Hilbert space $\mathcal{H}$ and $B$ is a boundedly
invertible linear operator on $\mathcal{H}$, the equality
\begin{equation}
\det(I+A)=\det(I+BAB^{-1})\label{eq.det.sim}%
\end{equation}
holds.

It can be shown that the determinant $\det(S_{r}(\lambda))$ extends to a
meromorphic function on the complex plane ($n$ odd) or the logarithmic plane
($n$ even) whose poles coincide, including multiplicity, with the resolvent
resonances. 

The real-valued function
\[
\sigma(\lambda)=\frac{1}{2\pi i}\log\det(S_{r}(\lambda))
\]
on $(0,\infty)$ is called the \emph{scattering phase} and behaves in many
respects analogously to the counting function for eigenvalues on a compact
manifold. For example, Christiansen \cite{Christiansen:1999} has shown that
the scattering phase for a class of scattering manifolds including those
considered here obeys the asymptotic law
\[
\sigma(\lambda)=-c_{n}\,\text{sc-}\operatorname*{vol}(X)\,\lambda
^{n}+\mathcal{O}(\lambda^{n-1})
\]
as $\lambda\rightarrow\infty$, where, in our case,
\[
\text{sc-}\operatorname*{vol}(X)=\lim_{\varepsilon\downarrow0}\left(
\int_{X_{\varepsilon}}dg-\frac{1}{n}\operatorname*{vol}(S^{n-1})\varepsilon
^{-n}-c_{n-1}\varepsilon^{n-1}-\cdots-c_{0}\log\varepsilon\right)
\]
The constant $c_{n}$ is the same constant that appears in Weyl's law for the
counting function of eigenvalues. The constants $c_{k}$ are chosen to make the
limit finite, and $X_{\varepsilon}$ is the compact set in $\mathbf{R}^{n}$
with $\left|  x\right|  \leq\varepsilon^{-1}$(equivalently,
sc-$\operatorname*{vol}(X)$ is the Hadamard finite part of
$\operatorname*{vol}_{g}(X_{\varepsilon})$ as $\varepsilon\downarrow0$). Note
that sc-$\operatorname*{vol}(X)$ may be positive, negative, or zero, depending
on $g$.  

\section{Technique for constructing isosphasal manifolds}

\label{sec.technique}

Before presenting the method we will use for constructing isosphasal metrics,
we review basic properties of group actions, in particular, torus actions.
Given an action of a compact Lie group $G$ on a manifold $M$, the
\emph{principal orbits} are the orbits with minimal isotropy. The union of the
principal orbits is an open dense subset $M^{\prime}$ of $M$. There exists a
subgroup $H$ of $G$ such the isotropy group of every element of $M^{\prime}$
is conjugate to $H$. Moreover, the isotropy group of an arbitrary element of
$M$ contains a subgroup conjugate to $H$. In case $G$ is a torus, it follows
that the isotropy group of every element contains $H$ itself. In particular,
if a torus action is effective, then $H$ is trivial and so the action on the
principal orbits is free. Thus $M^{\prime}$ is a principal $G$-bundle.

\begin{notation}
\label{notation.t} Suppose a torus $T$ acts smoothly on a connected manifold
$M$. For each character $\alpha:T\to S^{1}$ (where $S^{1}$ is the unit circle
in $\mathbf{C}$), write
\[
\mathcal{H}_{\alpha}=\{f\in C^{\infty}(M):f(z\cdot x) =\alpha(z)f(x)\text{ for
all }x\in M, z\in T\}.
\]
For $K$ a subtorus of $T$ of codimension at most one, let $\mathcal{C}
^{\infty}(M)^{K}$ denote the space of $K$-invariant smooth functions on $M$.
Then
\[
\mathcal{C}^{\infty}(M)^{K}=\oplus_{K\subset\ker(\alpha)}\,\mathcal{H}%
_{\alpha}.
\]
Thus by Fourier decomposition, we may decompose $\mathcal{C}^{\infty}(M)$ as
\[
\mathcal{C}^{\infty}(M)=\mathcal{C}^{\infty}(M)^{T}\oplus\left(  \oplus
_{K}\,\left(  \mathcal{C}^{\ \infty}(M)^{K}\ominus\mathcal{C}^{\infty}
(M)^{T}\right)  \right)
\]
where $K$ varies over all subtori of codimension one.
\end{notation}

\begin{proposition}
\label{prop.abstract method} \cite{G:2001}, \cite{Schueth:2002} Let $T$ be a
torus. Suppose $T$ acts effectively by isometries on two connected Riemannian
manifolds $(M_{1},g_{1})$ and $(M_{2},g_{2})$ and that the action of $T$ on
the principal orbits is free. Let $\mathring{M}^{\prime}_{i}$ be the union of
all principal orbits in $M_{i}$, so $\mathring{M}_{i}$ is an open, dense
submanifold of $M_{i}$ and a principal $T$-bundle, $i=1,2$. For each subtorus
$K$ of $T$ of codimension at most one, suppose that there exists a
$T$-equivariant volume-preserving diffeomorphism $F_{K}:M_{1}\rightarrow
M_{2}$ that induces an isometry $\bar{F}_{K}$ between the induced metrics on
the quotient manifolds $K\backslash\mathring{M}^{\prime}_{1}$ and
$K\backslash\mathring{M}^{\prime}_{2}$. With respect to the Fourier
decompositions of $\mathcal{C}^{\infty}(M_{1})$ and $\mathcal{C}^{\infty
}(M_{2})$ given in Notation \ref{notation.t}, let
\[
Q=F_{T}^{\ast}\oplus(\oplus_{K}F_{K}^{\ast}).
\]
Then
\[
\Delta_{1}=Q\circ\Delta_{2}\circ Q^{-1}
\]
where $\Delta_{i}$ is the Laplace operator on $\mathcal{C}^{\infty}(M_{i})$.
\end{proposition}

\begin{remark}
Proposition \ref{prop.abstract method} as stated here differs somewhat from
the original statements in \cite{G:2001} and \cite{Schueth:2002}. There the
manifolds $M_{i}$ were assumed to be compact and the conclusion was that they
were isospectral. However, the explicit intertwining operator $Q$ for the
Laplacians, which was constructed in the proof, did not use the assumption
that the manifolds were compact. A second difference between the statement
here and that in \cite{G:2001} is that a hypothesis involving preservation by
the diffeomorphisms $F_{K}$ of the mean curvature of the fibers has been
replaced by the condition that these diffeomorphisms be volume-preserving.
Dorothee Schueth made this simplifying change in her version
\cite{Schueth:2002} of the proposition, observing that the former and latter
conditions are equivalent.
\end{remark}

\begin{theorem}
\label{thm.gen.method}Suppose a torus $T$ acts effectively on $\mathbf{R}^{n}$
by orthogonal transformations. Let $g_{1}$ and $g_{2}$ be compact
perturbations of the Euclidean metric on $\mathbf{R}^{n}$ invariant under the
action of $T$. Assume that the Riemannian measure defined by both metrics
coincides with Lebesgue measure. Let $(\mathbf{\mathring{R}}^{n})^{\prime}$ be
the union of all principal orbits of the torus action. For each subtorus $K$
of $T$ of codimension at most one, suppose that there exists an orthogonal
transformation $F_{K}\in\mathrm{{O}(\mathbf{R}^{n})}$ commuting with $T$ which
induces an isometry $\overline{F_{K}}$ between the metrics induced by $g_{1}$
and $g_{2}$ on the quotient manifold $K\backslash(\mathbf{\mathring{R}}%
^{n})^{\prime}$. Then $g_{1}$ and $g_{2}$ have the same scattering phase.
\end{theorem}

\begin{remark}
We have stated the theorem only in the form needed for the examples given
here. However, the theorem may be generalized to other settings.
\end{remark}

\begin{proof}
The manifolds $(\mathbf{R}^{n},g_{1})$ and $(\mathbf{R}^{n},g_{2})$ satisfy
the hypotheses of Proposition \ref{prop.abstract method}. Define $Q:C^{\infty
}(\mathbf{R}^{n})\rightarrow C^{\infty}(\mathbf{R}^{n})$ as in the conclusion
of the proposition so that $\Delta_{1}=Q\circ\Delta_{2}\circ Q^{-1}$. Since
the $F_{K}$ are orthogonal maps of $\mathbf{R}^{n}$, the map $Q$ induces a map
$Q_{\partial}:C^{\infty}(S^{n-1})\rightarrow C^{\infty}(S^{n-1})$ which
extends to an invertible isometry of $L^{2}(S^{n-1})$. Since all orthogonal
maps commute with the antipodal map of $S^{n-1}$, the map $Q_{\partial}$
commutes with $S_{0}(\lambda)$ as defined in equation 2.2.

From its construction, it is clear that the intertwining map $Q$ preserves the
form of asymptotic expansions (\ref{eq.scatt.asy}). Moreover, if $u$ is a
solution $(\Delta_{g_{1}}-\lambda^{2})u=0$ having an asymptotic expansion of
the form (\ref{eq.scatt.asy}), then $Qu$ is a solution of $(\Delta_{g_{2}%
}-\lambda^{2})v=0$ having an asymptotic expansion of the form
\[
v(r\omega)=r^{(1-n)/2}e^{i\lambda r}\,h_{+}(\omega)+r^{(1-n)/2}e^{-i\lambda
r}h_{-}(\omega)+\mathcal{O}(r^{-(n+1)/2})
\]
as $r\rightarrow\infty$, where
\[
h_{+}(\omega)=(Q_{\partial}f_{+})(\omega)
\]
and
\[
h_{-}(\omega)=(Q_{\partial}f_{-})\left(  \omega\right)  .
\]
Let $S_{g_{1}}(\lambda)$ and $S_{g_{2}}(\lambda)$ be the scattering matrices
associated, respectively, to $(\mathbf{R}^{n},g_{1})$ and $(\mathbf{R}%
^{n},g_{2})$. Since
\[
f_{+}=S_{g_{1}}(\lambda)f_{-}%
\]
and
\[
h_{+}=S_{g_{2}}(\lambda)h_{-},
\]
it follows from the uniqueness statement in Proposition \ref{prop.unique}
that
\[
Q_{\partial}S_{g_{1}}(\lambda)f_{-}=S_{g_{2}}(\lambda)Q_{\partial}f_{-}.
\]
Since this holds for any $f_{-}\in\mathcal{C}^{\infty}(S^{n-1})$, and
$Q_{\partial}$ is an invertible linear map, we have
\[
S_{g_{2}}(\lambda)=Q_{\partial}S_{g_{1}}(\lambda)Q_{\partial}^{-1}.
\]
Since $Q_{\partial}$ and $Q_{\partial}^{-1}$ commute with the operator
$S_{0}(\lambda)$, we conclude that
\[
Q_{\partial}S_{g_{1}}(\lambda)S_{0}(\lambda)^{-1}Q_{\partial}^{-1}=S_{g_{2}%
}(\lambda)S_{0}(\lambda)^{-1},
\]
so that, on taking logarithms of determinants and using (\ref{eq.det.sim}),
\[
\sigma_{g_{2}}(\lambda)=\sigma_{g_{1}}(\lambda).
\]
\end{proof}

\section{Examples}

\label{sec.examples}

In \cite{G:2001}, the first author constructed continuous families of
Riemannian metrics on $\mathbf{R}^{n}$, $n\geq9$, which pairwise satisfy the
hypotheses of Theorem \ref{thm.gen.method} modulo the condition that the
metrics be Euclidean outside of a compact set. Dorothee Schueth
\cite{Schueth:2002} pointed out that the metrics could be modified to satisfy
this additional condition; in fact they could be flat outside of a compact set
of arbitrarily small volume. Moreover, Schueth constructed new continuous
families of metrics on $\mathbf{R}^{n}$, $n\geq8$, pairwise satisfying the
condition of Theorem \ref{thm.gen.method}. (Note that she lowered the minimum
dimension by one.) Additionally, Schueth constructed pairs, though not
continuous families, of such metrics on $\mathbf{R}^{6}$. In both these
papers, the focus was on compact manifolds. The metrics, once constructed,
were restricted to the unit ball and sphere. Using Proposition
\ref{prop.abstract method}, these restricted metrics were seen to be
isospectral. In the present context, we will conclude from Theorem
\ref{thm.gen.method} that the families of metrics on $\mathbf{R}^{n}$
constructed in these two papers are isosphasal.

We now review the construction of the metrics in \cite{G:2001} modified as in
\cite{Schueth:2002}.

\begin{definition}
\label{brack} {(i)} We will say that two skew-symmetric bilinear maps $\left[
\,\cdot\,,\,\cdot\,\right]  $ and $\left[  \,\cdot\,,\,\cdot\,\right]
^{\prime}$ taking $\mathbf{R}^{m}\times\mathbf{R}^{m}$ to $\mathbf{R}^{k}$ are
\emph{isospectral} if for each $Z\in\mathbf{R}^{k}$ there is an orthogonal
transformation $A_{Z}$ with the property that for every pair of vectors
$(x,y)\in\mathbf{R}^{m}\times\mathbf{R}^{m}$,
\[
\left\langle \left[  x,y\right]  ^{\prime},Z\right\rangle _{\mathbf{R}^{k}
}=\left\langle \left[  A_{Z}x,A_{Z}y\right]  ,Z\right\rangle _{\mathbf{R}^{k}
},
\]
where, here and in what follows, $\left\langle \,\cdot\,,\,\cdot
\,\right\rangle _{\mathbf{R}^{k}}$ denotes the Euclidean inner product on
$\mathbf{R}^{k}$. \newline {(ii)} We will say that the skew-symmetric bilinear
maps $\left[  \,\cdot\,,\,\cdot\,\right]  $ and $\left[  \,\cdot
\,,\,\cdot\,\right]  ^{\prime}$ are \emph{equivalent} if there exists an
orthogonal transformation $A$ of $\mathbf{R}^{m}$ and an orthogonal
transformation $C$ of $\mathbf{R}^{k}$, which preserves the lattice
$(2\pi\mathbf{Z})^{k}$, such that $\left\langle \left[  Ax,Ay\right]
^{\prime},Z\right\rangle =\left\langle \left[  x,y\right]  ,CZ\right\rangle $
for all $x,y\in\mathbf{R}^{m}$ and $Z\in\mathbf{R}^{k}$. We will also say that
the pair of maps $(A,C)$ is an equivalence of $\left[  \,\cdot\,,\,\cdot
\,\right]  $ and $\left[  \,\cdot\,,\,\cdot\,\right]  ^{\prime}$ in this case.
\end{definition}

\begin{remark}
\label{jmaps} {(i)} Our notation differs from that of \cite{G:2001}. The
bilinear maps $\left[  \,\,\cdot\,,\,\cdot\,\,\right]  :\mathbf{R}^{m}
\times\mathbf{R}^{m}\rightarrow\mathbf{R}^{k}$ correspond to linear maps
$j:\mathbf{R}^{k}\rightarrow\mathfrak{so}(m)$ via
\[
\left\langle \lbrack x,y],Z\right\rangle =\left\langle j(Z)x,y\right\rangle
\]
for all $x,y\in\mathbf{R}^{m}$ and $Z\in\mathbf{R}^{k}$. We say that $j$ and
$j^{\prime}$ are \emph{isospectral} if $j^{\prime}(Z)$ and $j(Z)$ are
isospectral linear operators for each $z\in\mathbf{R}^{k}$. We will say that
$j$ and $j^{\prime}$ are \emph{equivalent} if there exist orthogonal maps $A$
of $\mathbf{R}^{m}$ and $C$ of $\mathbf{R}^{m}$ such that $C$ preserves the
lattice $(2\pi\mathbf{Z})^{k}$ and such that $Aj^{\prime}(Z)A^{-1} =j(CZ)$ for
all $z\in\mathbf{R}^{k}$. These conditions correspond to the isospectrality
and equivalence conditions in Definition \ref{brack}. The $j$ maps were used
in \cite{G:2001} rather than the bracket maps $\lbrack\,,\,\rbrack$.
\newline {(ii)} Our notion of equivalence differs slightly from that in
\cite{G:2001} and \cite{GW:1997} in that we require $C$ to preserve the
lattice $(2\pi\mathbf{Z})^{k}$. This condition is added so that $C$ induces a
transformation of the torus $(2\pi\mathbf{Z})^{k} \backslash\mathbf{R}^{k}$.
\end{remark}

\begin{definition}
\label{metric} {(i)} Let $T$ be the $k$-torus $(2\pi\mathbf{Z})^{k}%
\backslash\mathbf{R}^{k}$ embedded in SO$(2k)$ as SO$(2)\times\dots\times
$SO$(2)$. Then $T$ acts on $\mathbf{R}^{2k}=\mathbf{R}^{2}\times\dots
\times\mathbf{R}^{2}$ by the standard SO$(2)$-action in each factor. This
action is not free but is inner-product preserving. The Lie algebra of
$T\,$\ is $\mathfrak{z} =\mathfrak{so}(2)\oplus\dots\oplus\mathfrak{so}
(2)\simeq\mathbf{R}^{k}$. Given $Z\in\mathfrak{z}$, define a vector field
$Z^{\ast}$ on $\mathbf{R}^{2k}$ by
\begin{equation}
Z_{u}^{\ast}=\left.  \frac{d}{dt}\right|  _{t=0}\left(  \exp(tZ)\cdot
u\right)  . \label{eq.z.field}%
\end{equation}
for $u\in\mathbf{R}^{2k}$. Observe that, for $x,y\in\mathbf{R}^{m}$,
$[x,y]\in\mathbf{R}^{k}=\mathfrak{z}$ so that we may use (\ref{eq.z.field}) to
define a vector field $[x,y]^{\ast}$ on $\mathbf{R}^{2k}$. \newline {(ii)}
Given a bilinear map $\left[  \,\cdot\,,\,\cdot\,\right]  :\mathbf{R}
^{m}\times\mathbf{R}^{m}\rightarrow\mathbf{R}^{k}$ and a smooth, compactly
supported function $\varphi:[0,\infty)\times\lbrack0,\infty)\rightarrow
\lbrack0,\infty)$, we now construct a Riemannian metric $g =g^{\left[
\,\cdot\,,\,\cdot\,\right]  ,\varphi}$ on $\mathbf{R}^{m+2k}$. Denote elements
of $\mathbf{R}^{m+2k}$ by $(x,u)$ with $x\in\mathbf{R}^{m}$ and $u\in
\mathbf{R}^{2k}$. First define $\psi:\mathbf{R}^{m+2k}\to\lbrack0,\infty)$ by
$\psi(x,u)=\varphi(\left\|  x\right\|  ^{2},\left\|  u\right\|  ^{2})$. For
$(x,u)\in\mathbf{R}^{m+2k}$, denote by $(Y,W)$ a typical element of the
tangent space $T_{(x,u)}\mathbf{R}^{m+2k}$, where, by standard
identifications, $Y\in\mathbf{R}^{m}$ and $W\in\mathbf{R}^{2k}$. We set
\[
g((0,W),(0,V))=\left\langle W,V\right\rangle _{\mathbf{R}^{2k}}
\]
and define the $g$-orthogonal complement to $\left\{  0\right\}
\oplus\mathbf{R}^{2k}$ in $T_{(x,u)}\mathbf{R}^{m+2k}$ as follows. For
$(x,u)\in\mathbf{R}^{m}\times\mathbf{R}^{2k}$ and $Y\in T_{x}\mathbf{R}^{m}$,
let
\[
\tilde{Y}_{x,u}=(Y,Z)
\]
with
\[
Z=\psi(x,u)\left[  x,Y\right]  _{u}^{\ast}.
\]

The $g$-orthogonal complement to $\left\{  0\right\}  \oplus\mathbf{R}^{2k}$
is taken to be
\[
\left\{  \tilde{Y}_{x,u}:Y\in T_{x}\mathbf{R}^{m}\right\}  .
\]
We put an inner product on this space so that the map $Y\mapsto\tilde{Y}
_{x,u}$ is an isometry where $\mathbf{R}^{m}$ has the Euclidean inner product.
\end{definition}

Note that, for $(x,u)$ outside of the support of $\psi$, we have $\tilde
{Y}_{x,u}=Y$. Thus the metric so constructed is identical to the Euclidean
metric away from the support of $\psi$.

\begin{proposition}
\label{prop.isosp} Suppose that $\left[  \,\cdot\,,\,\cdot\,\right]  $ and
$\left[  \,\cdot\,,\,\cdot\,\right]  ^{\prime}$ are isospectral in the sense
of Definition \ref{brack} and that $g$ and $g^{\prime}$ are metrics
constructed as in Definition \ref{metric} from the data $(\varphi,\left[
\,\cdot\,,\,\cdot\,\right]  )$ and $(\varphi,\left[  \,\cdot\,,\,\cdot
\,\right]  ^{\prime})$ for the \emph{same} nonnegative function $\varphi
\in\mathcal{C}_{0}^{\infty}(\mathbf{R}^{+}\times\mathbf{R}^{+})$. Then $g$ and
$g^{\prime}$ are isophasal.
\end{proposition}

\begin{proof}
We apply Theorem \ref{thm.gen.method}. Let $W$ denote the union of the
principal orbits for the action of $T$ on $\mathbf{R}^{2k}$. By identifying
$\mathbf{R}^{2k}$ with $\mathbf{C}^{k}$, we may write
\[
W=\left\{  (z_{1},\dots,z_{k})\in\mathbf{C}^{k}:z_{1}\neq0,\dots,z_{k}%
\neq0\right\}  .
\]
The union of the principal orbits for the action of $T$ on $\mathbf{R}^{m+2k}$
is given by $\mathbf{R}^{m}\times W$.

If $K\subset T$ is a subtorus of codimension one, then in the Lie algebra
$\mathfrak{z}$ of $T$, there is a vector $Z$ orthogonal to the Lie subalgebra
$\mathfrak{k}$ of $K$. By hypothesis, there is an orthogonal transformation
$A_{Z}\in$O$(m)$ so that
\[
\left\langle \left[  x,y\right]  ^{\prime},Z\right\rangle _{\mathbf{R}^{k}
}=\left\langle \left[  A_{Z}x,A_{Z}y\right]  ,Z\right\rangle _{\mathbf{R}%
^{k}}
\]
for any $x$ and $y$ belonging to $\mathbf{R}^{m}$. Letting $g_{K}$ and
$g_{K}^{\prime}$ be the metrics on $K\backslash(\mathbf{\mathring{R}}%
^{m}\times W)$ induced by $g$ and $g^{\prime}$, it follows from Definition
\ref{metric}(ii) that the orthogonal map
\[
\tau_{K}(x,u)=(A_{Z}x,u)
\]
of $\mathbf{R}^{m+2k}$ induces an isometry from $(K\backslash
(\mathbf{\mathring{R}}^{m}\times W),g_{K})$ to $(K\backslash(\mathbf{\mathring
{R}}^{m}\times W),g_{K}^{\prime})$. Thus the hypotheses of Theorem
\ref{thm.gen.method} are satisfied, and we conclude that the metrics are isophasal.
\end{proof}

\begin{remark}
In \cite{G:2001}, the function $\psi$ did not appear; \ i.e., $\varphi$ (and
thus $\psi$) was identically one. As mentioned above, it was Dorothee Schueth
that realized the function $\psi$ could be inserted so that the metrics are
Euclidean outside of a compact set.
\end{remark}

By referring to the proofs of Proposition \ref{prop.isosp} and of Theorem
\ref{thm.gen.method}, we can give an explicit description of the intertwining
operators between the Laplacians of the metrics in Proposition
\ref{prop.isosp} and between their scattering phases as follows:

\begin{proposition}
\label{prop.intertwine} Define $g$ and $g^{\prime}$ as in Proposition
\ref{prop.isosp}. Writing $\mathbf{R}^{2k}=\mathbf{R}^{2}\times\dots
\times\mathbf{R}^{2}$ and letting $(r_{i},\theta_{i})$ denote polar
coordinates on the $i$th factor, we obtain coordinates $(x,r,\theta)$ on
$\mathbf{R}^{m+2k}$, where $x =(x_{1},\dots,x_{m})$, $r=(r_{1},\dots,r_{k})$
and $\theta=(\theta_{1},\dots,\theta_{k})$. For $Z\in\mathbf{R}^{k}$, choose
$A_{Z}$ as in Definition \ref{brack}. Define $Q:C^{\infty}(\mathbf{R}%
^{m+2k})\to C^{\infty}(\mathbf{R}^{m+2k})$ by
\[
Q(f)(x,r,\theta)=\Sigma_{Z\in(\frac{1}{2\pi}\mathbf{Z})^{n}}(\int
_{[0,2\pi]^{n}}\,f(A_{Z}(x),r,\sigma)e ^{-iZ\cdot\sigma}d\sigma)\,e^{iZ\cdot
\theta}.
\]
Then $Q$ intertwines the Laplacians of the metrics $g$ and $g^{\prime}$ on
$\mathbf{R}^{m+2k}$, and the associated map $Q_{\partial}$, defined as in
Theorem \ref{thm.gen.method}, intertwines their scattering phases.
\end{proposition}

We now consider whether these metrics are isometric.

\begin{proposition}
\label{prop.bundle.isom} Fix $\varphi$. Let $g$ and $g^{\prime}$ be the
metrics defined as in Definition \ref{metric} by nontrivial maps
$[\,\cdot\,\,,\,\cdot\,\,]$ and $[\,\,\cdot\,,\,\cdot\,\,]^{\prime}%
:\mathbf{R}^{m}\times\mathbf{R}^{m}\rightarrow\mathbf{R}^{k}$ together with
$\varphi$.

\item {(i)} Suppose that $\tau$ is an isometry from $(\mathbf{R}^{m+2k},g)$ to
$(\mathbf{R}^{m+2k},g^{\prime})$ that carries $T$-orbits to $T$-orbits. Then
$\tau$ is of the form
\[
(\tau(x,u)=(A(x),\tilde{C}(u))
\]
where $A\in$O$(m)$, $\tilde{C}\in$O$(2k)$, and $\tilde{C}$ normalizes $T$.
Letting $C$ be the automorphism of the Lie algebra $\mathbf{R}^{k}
=\mathfrak{so}(2)\oplus\cdots\oplus\mathfrak{so}(2)$ of $T$ given by
conjugation by $\tilde{C}$, then the pair $(A,C)$ is an equivalence of
$[\,\cdot\,\,,\,\cdot\,\,]$ and $[\,\cdot\,\,,\,\cdot\,\,]^{\prime}$, as in
Definition 4.1.

\item {(ii)} Conversely, every map $\tau$ of this form is an isometry between
the two metrics.
\end{proposition}

\begin{proof}
(ii) is straightforward and is left to the reader.

(i) Since the isometry $\tau$ carries $T$-orbits to $T$-orbits, it must
preserve the open dense subset $\mathbf{R}^{m}\times W$, where $W$ is given in
the proof of Proposition \ref{prop.isosp}. The submanifold $\mathbf{R}%
^{m}\times W$ has the structure of a principal $T$ bundle over $\mathbf{R}%
^{m}\times(T\backslash W)\simeq\mathbf{R}^{m} \times(\mathbf{R}^{+})^{k}$. The
metrics $g$ and $g^{\prime}$ both induce the standard Euclidean metric on the
quotient $\mathbf{R}^{m}\times(\mathbf{R} ^{+})^{k}$. The isometry $\tau$
induces an isometry $\bar{\tau}$ of $\mathbf{R}^{m}\times(\mathbf{R}^{+})^{k}%
$. Such an isometry is the composition of a translation in $\mathbf{R}^{m}$
with an orthogonal transformation of the form $A\times P$, where $A\in
\text{O}(m)$ and $P$ permutes the coordinates in $(\mathbf{R}^{+})^{k}$. We
claim that the translation factor is trivial. To see this, note that the
metrics $g$ and $g^{\prime}$ on $\mathbf{R}^{m+2k}$ are Euclidean on the
complement of $\{(x,u)\in\mathbf{R}^{m+2k}:(\Vert x\Vert,\Vert u\Vert
)\in\text{supp}(\varphi)\}$. Letting $R$ be minimal such that $\text{supp}%
(\varphi)\subset\{(s,t):s^{2}+t^{2}\leq R\}$, then $g$ and $g^{\prime}$ are
Euclidean on the region $\{(x,u):\Vert x\Vert^{2}+\Vert u\Vert^{2}\geq R\}$
and not on any translate of this region. Hence $\bar{\tau}$ must preserve the
image of this region in $\mathbf{R}^{m}\times(\mathbf{R}^{+})^{k}$, and the
claim follows.

For each $x\in\mathbf{R}^{m}$, $\tau$ restricts to an isometry from the
Euclidean space $\{x\}\times\mathbf{R}^{2k}$ to the Euclidean space
$\{A(x)\}\times\mathbf{R}^{2k}$. Canonically identifying both spaces with
$\mathbf{R}^{2k}$, this isometry preserves the origin, since the origin is the
unique $T$-orbit which is a single point. Thus $\tau$ is of the form
$\tau(x,u)=(A(x),B_{x}(u))$ with $B_{x}\in\text{O}(2k)$ for each
$x\in\mathbf{R}^{m}$. We may identify $T$ with the maximal torus
$T=\text{SO}(2)\times\dots\times\text{SO}(2)$ of $\text{O}(2k)$. Since $\tau$
carries $T$-orbits to $T$-orbits, each $B_{x}$ must normalize $T$. Noting that
$T$ has finite index in its normalizer in $\text{O}(2k)$ and that $B_{x}$
depends smoothly on $x$, there must exist $\tilde{C}\in\text{O}(2k)$,
independent of $x$, and $z(x)\in T$ such that $B_{x}=z(x)\circ\tilde{C}$. The
permutation $P$ in the expression for $\bar{\tau}$ is the map of $T\backslash
W$ induced by $\tilde{C}$.

We next show that $(A,C)$ defines an equivalence of $\left[  \,\cdot
\,,\,\cdot\,\right]  $ and $\left[  \,\cdot\,,\,\cdot\,\right]  ^{\prime}$,
where $C$ is defined from $\tilde{C}$ as in the statement of the proposition.
Observe that $\tau\circ\psi=\tau$, since $\tau(x,u)$ preserves the norms of
the two coordinates. Since $\tau_{\ast}$ maps $g$-horizontal vectors at each
point $(x,u)$ (i.e., vectors $g$-orthogonal to the orbit of $T$ through
$(x,u)$) to $g^{\prime}$-horizontal vectors through $\tau(x,u)$, we have for
$y\in\mathbf{R}^{m}$, viewed as a tangent vector to $\mathbf{R}^{m}$ at $x$,
\[
\tau_{\ast(x,u)}(y+\psi(x,u)[x,y]_{(x,u)}^{\ast})=A(y)+\psi
(x,u)[A(x),A(y)]_{(\tau(x,u)}^{^{\prime}\ast}.
\]
Noting that $\tau_{\ast}(Z^{\ast})=C(Z^{\ast})$ for $Z$ in the Lie algebra
$\mathfrak{z}=\mathbf{R}^{k}$ of $T$, it follows that $[A(x),A(y)]^{^{\prime
}\ast}=(C([x,y])^{\ast}$. Since the map $Z\rightarrow Z^{\ast}$ is injective
on $\mathfrak{z}$, we see that $(A,C)$ is an equivalence of $\left[
\,\cdot\,,\,\cdot\,\right]  $ and $\left[  \,\cdot\,,\,\cdot\,\right]
^{\prime}$.

It remains to show that the map $x\rightarrow z(x)$ is constant. Fix a point
$x_{0}$ and let $z=z(x_{0})$. Define $\mu(x,u) =(A(x),z\cdot\tilde{C}(u))$. By
(ii), $\mu$ is an isometry from $g$ to $g^{\prime}$. Hence $\alpha: =\tau
^{-1}\circ\mu$ is an isometry of $g$ of the form $(x,u)\rightarrow(x,w(x)\cdot
u)$ for some map $w:\mathbf{R}^{m}\rightarrow T$ satisfying $w(x_{0})=1$,
where $1$ denotes the identity element in $T$. At points of the form
$(x_{0},u))$, $u\in\mathbf{R}^{2k}$, the differential $\alpha_{\ast}$ acts as
the identity both on the tangent space to the fiber and on the horizontal
space. Thus $\alpha_{\ast(x_{0},u)}=Id$. Since an isometry is uniquely
determined by its value and its differential at a single point, it follows
that $\alpha=Id$, i.e., that $z(x)\equiv z$. Replacing $\tilde{C}$ by
$z\cdot\tilde{C}$. The lemma follows.
\end{proof}

\begin{corollary}
\label{cor.com.isom} Let $g$ be the Riemannian metric on $\mathbf{R}^{m+2k}$
defined from the data $(\left[  \,\cdot\,,\,\cdot\,\right]  ,\varphi)$ as in
Definition \ref{metric}. Then the centralizer of $T$ in the group of all
isometries of $g$ consists of all maps $\tau$ of $\mathbf{R}^{m+2k}$ of the
form $\tau(x,u)=(A(x),z\cdot u)$ such that $A\in\text{O}(m)$ preserves
$\left[  \,\cdot\,,\,\cdot\,\right]  $ (i.e., the pair $(A,\text{Id})$ is a
self-equivalence of $\left[  \,\cdot\,,\,\cdot\,\right]  $) and such that
$z\in T$.
\end{corollary}

\begin{proof}
An isometry that commutes with $T$ must carry $T$-orbits to $T$-orbits. Thus
the corollary follows from Proposition \ref{prop.bundle.isom} and the fact
that $T$ is its own centralizer in O$(2k)$. .
\end{proof}

\begin{proposition}
\label{prop.isom} Suppose that $g$ and $g^{\prime}$ are metrics constructed
from the data $\left(  \left[  \,\cdot\,,\,\cdot\,\right]  ,\varphi\right)  $
and $\left(  \left[  \,\cdot\,,\,\cdot\,\right]  ^{\prime},\varphi\right)  $
as in Definition \ref{metric}, where $\left[  \,\cdot\,,\,\cdot\,\right]  $
and $\left[  \,\cdot\,,\,\cdot\,\right]  ^{\prime}$ are inequivalent in the
sense of Definition \ref{brack}. Assume that $\left[  \,\cdot\,,\,\cdot
\,\right]  $ satisfies the genericity condition that $\left[  \,\cdot
\,,\,\cdot\,\right]  $ is invariant under only finitely many orthogonal
transformations of $\mathbf{R}^{m}$. Then $g$ is not isometric to $g^{\prime}$.
\end{proposition}

\begin{proof}
By Corollary \ref{cor.com.isom} and the genericity condition on $\left[
\,\cdot\,,\,\cdot\,\right]  $, $T$ is a maximal torus in the full isometry
group of $g$. Now suppose that $\rho:(\mathbf{R}^{n},g)\rightarrow
(\mathbf{R}^{n},g^{\prime})$ is an isometry. Since the metrics are isometric,
$T$ must also be a maximal torus in the full isometry group of $g^{\prime}$.
By the conjugacy of the maximal tori in any Lie group, we may assume after
composing with an isometry of $g^{\prime}$ that $\rho$ carries $T$-orbits to
$T$-orbits. By Proposition \ref{prop.bundle.isom}, it follows that $\left[
\,\cdot\,,\,\cdot\,\right]  $ and $\left[  \,\cdot\,,\,\cdot\,\right]
^{\prime}$ are equivalent, contradicting the hypothesis.
\end{proof}

By Proposition \ref{prop.isosp} and Proposition \ref{prop.isom}, if $\left[
\,\cdot\,,\,\cdot\,\right]  $ and $\left[  \,\cdot\,,\,\cdot\,\right]
^{\prime}$ are isospectral, inequivalent maps $\mathbf{R}^{m}\times
\mathbf{R}^{m}\rightarrow\mathbf{R}^{k}$ in the sense of Definition 4.1 and if
$\left[  \,\cdot\,,\,\cdot\,\right]  $ satisfies the genericity condition of
Proposition \ref{prop.isom}, then the metrics on $\mathbf{R}^{m+2k}$
constructed from the data $(\left[  \,\cdot\,,\,\cdot\,\right]  ,\varphi)$ and
$(\left[  \,\cdot\,,\,\cdot\,\right]  ^{\prime},\varphi)$ as in Definition
\ref{metric}, for any fixed choice of $\varphi$, are isospectral but not
isometric. The following lemma shows that such isospectral, inequivalent maps
are plentiful.

\begin{proposition}
\cite{GW:1997}\label{prop.exist} Let $k=2$, and let $m$ be any positive
integer other than $1,2,3,4$, or $6$. Let $W_{m}$ be the real vector space
consisting of all anti-symmetric bilinear maps from $\mathbf{R}^{m}
\times\mathbf{R}^{m}$ to $\mathbf{R}^{k}$. Then there is a Zariski open subset
${\mathcal{O}}_{m}$ of $W_{m}$ (i.e., ${\mathcal{O}}_{m}$ is the complement of
the zero locus of some non-zero polynomial function on $W$) such that each
$\left[  \,\cdot\,,\,\cdot\,\right]  \in{\mathcal{O}}_{m}$ belongs to a
$d$-parameter family of isospectral, inequivalent elements of $W_{m}$. Here
$d\geq\frac{m(m-1)}{2}-[\frac{m}{2}]([\frac{m}{2}]+2)>1$. In particular, $d$
is of order at least $\mathcal{O}(m^{2})$. Moreover, the elements of
${\mathcal{O}}_{m}$ satisfy the genericity condition of Proposition
\ref{prop.isom}.
\end{proposition}

The statement of the proposition in \cite{GW:1997} is in the language of
Remark \ref{jmaps}. The final statement of the proposition was not explicitly
stated in \cite{GW:1997}; however, a glance at the proof given there shows
that the genericity condition is one of the defining properties of the Zariski
open set ${\mathcal{O}}_{m}$ constructed there. While the proposition omits
$m=6$, an explicit example of a continuous family of isospectral, inequivalent
maps $\mathbf{R}^{m}\times\mathbf{R}^{m}\rightarrow\mathbf{R}^{2}$ was also
constructed in \cite{GW:1997}.

We have now proven Theorem \ref{thm.main} when $n\geq9$. (For $n=10$, we refer
to the comment immediately above.)

To prove Theorem \ref{thm.main} in the case $n=8$ for continuous families and
$n=6$ for pairs, we refer to the article \cite{Schueth:2002} by Dorothee
Schueth. There Schueth constructed metrics on $\mathbf{R}^{n}$ from data
$(L,\psi)$ consisting of a particular type of linear map $L$ and a cut-off
function $\psi$ on $\mathbf{R}^{n}$ of the same type used above. (The maps
$L$, which play the role of the $\left[  \,\cdot\,,\,\cdot\,\right]  $ maps in
the construction above, are denoted $j$ or $c$ in the two different
constructions given in \cite{Schueth:2002}.) She defined notions of
isospectrality and equivalence of the linear maps. An argument analogous to
the proof of Proposition \ref{prop.isosp} shows that, for fixed $\psi$, the
metrics on $\mathbf{R}^{n}$ constructed from isospectral linear maps $L$ and
$L^{\prime}$ are isosphasal. To discuss the condition for non-isometry, we
will for simplicity require that the cut-off function $\psi$ be radial. The
metrics constructed in \cite{Schueth:2002} are Euclidean outside of the
support of $\psi$ but not on any open set on which $\psi$ is positive. Since
$\psi$ is supported on a ball $B$ about the origin, any isometry between the
metrics must therefore carry this ball to itself. Under a genericity condition
analogous to that in Proposition \ref{prop.isom}, Schueth proved that the
metrics on the ball are not isometric provided that the associated linear maps
are inequivalent. This completes the proof.

\begin{example}
\label{ex} We give an explicit triple of isophasal metrics on $\mathbf{R}%
^{12}$ and compare their geometries. We let $k=3$ and $m=6$. Define three maps
$\left[  \,\cdot\,,\,\cdot\,\right]  _{i}:\mathbf{R}^{6}\times\mathbf{R}^{6}
\rightarrow\mathbf{R}^{3}$ as follows.

To define $\left[  \,\cdot\,,\,\cdot\,\right]  _{1}$ and $\left[
\,\cdot\,,\,\cdot\,\right]  _{2}$, view $\mathbf{R}^{6}$ as $\mathbf{R}
^{3}\times\mathbf{R}^{3}$ and denote elements of $\mathbf{R}^{6}$ as ordered
pairs $(x,y)$, with $x,y\in\mathbf{R}^{3}$. Let $\times$ denote the cross
product on $\mathbf{R}^{3}$. Define
\[
\lbrack(x,y),(x^{\prime},y^{\prime})]_{1}=x\times x^{\prime}+y\times
y^{\prime}
\]
and
\[
\lbrack(x,y),(x^{\prime},y^{\prime})]_{2}=x\times x^{\prime}-y\times
y^{\prime}.
\]

To define $\left[  \,\cdot\,,\,\cdot\,\right]  _{3}$, view $\mathbf{R}^{6}$ as
$H\times\mathbf{R}^{2}$, where $H$ denotes the quaternions. Denote elements of
$\mathbf{R}^{6}$ as pairs $(q,y)$, with $q\in H$, $y\in\mathbf{R}^{2}$. View
the target space $\mathbf{R}^{3}$ as the purely imaginary quaternions. Define
\[
\lbrack(q,y),(q^{\prime},y^{\prime})]_{3}=\text{Im}(qq^{\prime}),
\]
where $qq^{\prime}$ is the quaternionic product.

To see that the three bracket maps are isospectral, it is easier to consider
the associated maps $j_{i}:\mathbf{R}^{3}\rightarrow\mathfrak{so}(6)$ defined
as in Remark \ref{jmaps}. We have
\[
j_{1}(z)(x,y)=(z\times x,z\times y),
\]
\[
j_{2}(z)(x,y)=(z\times x,-z\times y),
\]
and
\[
j_{3}(z)(q,y)=(zq,0),
\]
where in the final equation, $zq$ denotes quaternionic multiplication of the
purely imaginary quaternion $z$ with the quaternion $q$. In each case, the
eigenvalues of $j_{i}(z)$ are $\Vert z\Vert\sqrt{-1}$, $-\Vert z\Vert\sqrt
{-1}$, and $0$, each occurring with multiplicity 2. Thus $j_{1}(z)$,
$j_{2}(z)$ and $j_{3}(z)$ are similar transformations for each $z$, and hence
the $j_{i}$ are mutually isospectral. Equivalently, the $\left[
\,\cdot\,,\,\cdot\,\right]  _{i}$ are mutually isospectral. Thus fixing a
choice of $\varphi$, we obtain a triple of isosphasal metrics $g_{i}$ on
$\mathbf{R}^{12}$.

The isometry groups $\text{Iso}(g_{i})$ of the three metrics vary both in
their dimension and structure. By Corollary \ref{cor.com.isom} and Remark
\ref{jmaps}, every isometry of $g_{i}$ that commutes with $T$ is the
composition of an element of $T$ with an isometry of the form $A\times
\text{Id}$ acting on $\mathbf{R}^{12} =\mathbf{R}^{6}\times\mathbf{R}^{6}$,
where $A\in$O$(6)$ commutes with all the $j_{i}(z)$, $z\in\mathbf{R}^{k}$. The
image of $j_{1}$ in $\mathfrak{so}(6)$ is the set of all matrices of the form
\[
\left(
\begin{matrix}
B & 0\\
0 & B
\end{matrix}
\right)
\]
with $B\in\mathfrak{so}(3)$. The centralizer of this image in O$(6)$ is the
one-parameter subgroup (circle) generated by the skew-symmetric matrix
\[
\left(
\begin{matrix}
0 & \text{Id}\\
\text{-Id} & 0
\end{matrix}
\right)  .
\]
Thus the centralizer of $T$ in the connected component of $\text{Iso}(g_{1})$
is isomorphic to $T\times S^{1}$, a four-dimensional torus. In particular, the
maximal tori in $\text{Iso}(g_{1})$ are four-dimensional.

The image of $j_{2}$ in $\mathfrak{so}(6)$ is the set of all matrices of the
form
\[
\left(
\begin{matrix}
B & 0\\
0 & -B
\end{matrix}
\right)
\]
with $B\in\mathfrak{so}(3)$. The centralizer of this image in O$(6)$ is
trivial. Thus the three-dimensional torus $T$ is a maximal torus in
$\text{Iso}(g_{2})$ and is its own centralizer.

The connected component of the centralizer of the image of $j_{3}$ in
$\mathfrak{so}(6)$ is isomorphic to $\text{SU}(2)\times\text{SO}(2)$, where
the 3-sphere $\text{SU}(2)$ is identified with the unit quaternions acting on
$H$ by right multiplication and where $\text{SO}(2)$ acts on the
$\mathbf{R}^{2}$ factor. Since a maximal torus in $\text{SU}(2)\times
\text{SO}(2)$ is two-dimensional, the maximal tori in $\text{Iso}(g_{3})$ are
five-dimensional. Moreover, the semisimple group $\text{SU}(2)\times
\text{SO}(2)$ acts by isometries preserving the bundle structure.
\end{example}

\section{Existence of Resonances}

\label{sec.a2}

To show that our examples are not trivial, we will show that the metrics that
we have constructed have infinitely many scattering resonances, in contrast to
the Euclidean Laplacian on $\mathbf{R}^{n}$ which has none. We use the methods
of Sa Barreto-Tang \cite{SbT:1998} if $n$ is odd and Tang \cite{Tang:2000} if
$n$ is even.

In \cite{SbT:1998} it is shown that if the Laplacians of two metrics on
$\mathbf{R}^{n}$ ($n$ odd) which differ from the Euclidean metric by a
super-exponentially decaying perturbation have the same resonances, then they
also have the same heat invariants $a_{k}$ for $k\geq2$. This implies that a
metric with non-vanishing heat invariant $a_{2}$ must have at least finitely
many resonances, and a closer analysis of the renormalized wave trace shows
that, in fact, the number of resonances must be infinite in this case. For the
case of $n$ even, it is shown in \cite{Tang:2000} that, if a metric has only
finitely many resonances, then the heat invariants $a_{k}$ for $k\geq2$ must
vanish. These papers also give various geometric hypotheses under which
$a_{2}\neq0$.

Here, we carry out the proof that the heat invariant $a_{2}$ is non-zero for
the specific metrics constructed in Section 4. For convenience, we will
restrict to the cases that $k=2$ or $k=3$ in the notation of Section 4. These
two cases include all the examples constructed by Proposition \ref{prop.exist}
as well as Example \ref{ex}.

Fix a non-trivial bilinear map $\left[  \,\cdot\,,\,\cdot\,\right]
:{\mathbf{R}}^{m}\times\mathbf{R}^{m}\rightarrow\mathbf{R}^{k}$. For each
non-trivial compactly supported function $\varphi:[0,\infty)\times
\lbrack0,\infty)\rightarrow\lbrack0,\infty)$, define the metric $g^{\left[
\,\cdot\,,\,\cdot\,\right]  ,\varphi}$ on ${\mathbf{R}}^{m+2k}$ as in
\S\ref{sec.examples}. The goal of this section is to show that, for generic
$\varphi$, the scattering matrix for the metric $g^{\left[  \,\cdot
\,,\,\cdot\,\right]  ,\varphi}$ has infinitely many resonances. It suffices to
show that the second heat invariant $a_{2}(g^{\left[  \,\cdot\,,\,\cdot
\,\right]  ,\varphi})$ is non-zero.

Recall that for any Riemannian metric $g$ on an $n$-dimensional manifold $M$,
the second heat invariant is given by
\[
a_{2}(g) =\frac{(4\pi)^{-\frac{n}{2}}}{360}\int_{M}(5\tau^{2}-2\Vert
\operatorname{Ric}\Vert^{2}+2\Vert R\Vert^{2})\,\operatorname*{dvol}
\nolimits_{g},
\]
where $\tau$ denotes the scalar curvature.

\begin{notation}
\label{s} Given $\varphi: [0,\infty)\times[0,\infty)\to[0,\infty)$, consider
the family of functions $\varphi_{s}$, $s\in{\mathbf{R}}^{+}$, given by
$\varphi_{s}(t_{1},t_{2}) =\varphi(t_{1},s^{2}t_{2})$. Set
\[
g^{s}=g^{\left[  \,\cdot\,,\,\cdot\,\right]  ,\varphi_{s}}
\]
and denote by $R^{(s)}$ and $\operatorname{Ric}^{(s)}$ the curvature tensor
and Ricci tensor of the metric $g^{s}$.
\end{notation}

\begin{theorem}
\label{thm.5.1} Let $k=2$ or $3$, let $\left[  \,\cdot\,,\,\cdot\,\right]
:{\mathbf{R}}^{m}\times\mathbf{R}^{m}\rightarrow\mathbf{R}^{k}$ be a non-zero
bilinear map, and let $\varphi:[0,\infty)\times\lbrack0,\infty)\rightarrow
\lbrack0,\infty)$ be a non-trivial compactly supported $\mathcal{C}^{\infty}$
function. We use Notation \ref{s}. Except for possibly finitely many values of
$s$, $a_{2}(g^{s})\neq0$.
\end{theorem}

\begin{definition}
\label{def.5.3} Let $\{f^{s}\}$, $s\in{\mathbf{R}}^{+}$, be a one-parameter
family of functions on ${\mathbf{R}}^{m+2k}$. We will say that $\{f^{s}\}$ is
a {\textit{homogeneous}} $(r,s)$-{\textit{deformation of degree}} $d$ if
$f^{s}(x,r,\theta)=s^{d}f^{1}(x,{\tilde{r}},\theta)$ for all $s$, where
${\tilde{r}}=({\tilde{r}}_{1},\dots,{\tilde{r}}_{k}) =(sr_{1},\dots, sr_{k})$.
We will say $\{f^{s}\}$ is an $(r,s)$-{\textit{deformation of degree}} $d$ if
$f^{s}=f_{1}^{s}+\dots f_{l}^{s}$ where each $\{f_{i}^{s}\}$ is a homogeneous
$(r,s)$-deformation, say of degree $d_{i}$, and where $d=\max\{d_{i}\}$.
Assuming the $d_{i}$'s are distinct, we will refer to $f_{i}^{s}$ as the
homogeneous term of degree $d_{i}$ in $f^{s}$.
\end{definition}

The following two lemmas are elementary.

\begin{lemma}
\label{lemma.5.4} If $\{f^{s}\}$ is a homogeneous $(r,s)$-deformation of
degree $d$, then:\newline {(i)} $\{\frac{\partial}{\partial x_{j}}f^{s}\}$,
$1\leq j\leq m$, is also a homogeneous $(r,s)$-deformation of degree $d$ while
$\{\frac{\partial}{\partial r_{p}}f^{s}\}$, $1\leq p\leq k$, is a homogeneous
$(r,s)$-deformation of degree $d+1$;\newline {(ii)} $\{r_{p}f^{s}\}$, $1\leq
p\leq k$, is a homogeneous $(r,s)$-deformation of degree $d-1$.\newline
{(iii)} If $\{h^{s}\}$ is a homogeneous $(r,s)$-derivation of degree
$d^{\prime}$, then $\{f^{s}h^{s}\}$ is a homogeneous $(r,s)$-deformation of
degree $d+d^{\prime}$.
\end{lemma}

\begin{lemma}
\label{lemma.5.5} If $\{f^{s}\}$ is a homogeneous $(r,s)$-deformation of
degree $d$ and if each $f^{s}$ is continuous with compact support, then
\[
\int_{\mathbf{R}^{m+2k}}\,f^{s} =s^{d-2k}\int_{\mathbf{R}^{m+2k}}\,f^{1}
\]
where the integrals are with respect to Lebesgue measure.
\end{lemma}

To prove that the function $a_{2}(g^{s})$ is non-zero except for possibly
finitely many choices of $s$, we will show that the integrand is an
$(r,s)$-deformation of degree 2, and that the homogeneous term of degree two
is strictly positive. Consequently, after multiplying by an appropriate power
of $s$, the heat invariant $a_{2}(g^{s})$ is a non-trivial polynomial in $s$,
from which Theorem \ref{thm.5.1} follows.

\begin{notation}
\label{rem.5.2} (i) Write ${\mathbf{R}}^{2k}$ as ${\mathbf{R}}^{2}
\times\mathbf{R}^{2}\dots\times\mathbf{R}^{2}$ and let $(r_{i},\theta_{i})$,
$i=1,2, \dots,k$, denote the polar coordinates on the $k$ factors
${\mathbf{R}}^{2}$. We thus coordinatize a dense open subset of $\mathbf{R}%
^{m+2k}$ by $(x,r,\theta)=(x_{1},\dots,x_{m},r_{1},r_{2},\dots r_{k}%
,\theta_{1},\theta_{2},\dots, \theta_{k})$.

(ii) For $\{Z_{1}, \dots,Z_{k}\}$ the standard basis of $\mathfrak{z}%
=\mathbf{R}^{k}$, the vector field $Z_{p}^{\ast}$, $p=1,\dots,k$, defined in
equation (\ref{eq.z.field}), is given by $\frac{\partial}{\partial\theta_{p}}$.

(iii) Let indices $i,j,k$ run from $1$ to $m$, indices $p,q$ range over
$1,2,\dots,k$, and Greek indices range over $1,\dots,m+2k$. Define an
orthonormal frame field as follows: For $p=1,2,\dots, k$, let $\hat{r}%
_{p}=\frac{\partial}{\partial r_{p}}$ and let $\hat{\theta}_{p}$ be a unit
vector in the direction of $\frac{\partial}{\partial\theta_{p}}$, i.e.,
$\hat{\theta}_{p}=\frac{1}{r_{p}}\theta_{p}$. For $\{e_{1},\dots,e_{m}\}$ the
standard basis of ${\mathbf{R}}^{m}$, define
\[
a_{ip}^{s}(x,r)=\varphi_{s}(\Vert x\Vert^{2},\Vert r\Vert^{2})\langle\lbrack
x,e_{i}],Z_{p}\rangle.
\]
Set
\[
\hat{x}_{i}^{s} =e_{i}+a_{i1}^{s}(x,r)\frac{\partial}{\partial\theta_{1}
}+\dots+a_{ik}^{s}(x,r)\frac{\partial}{\partial\theta_{k}} =e_{i}+a_{i1}
^{s}(x,r)r_{1}\hat{\theta}_{1}+\dots+a_{ik}^{s}(x,r)r_{k}\hat{\theta}_{k}.
\]
Then
\[
\{E_{1}^{s},\dots,E_{m+2k}^{s}\} =\{\hat{x}_{1}^{s},\dots,\hat{x}_{m}^{s}
,\hat{r}_{1},\dots,\hat{r}_{k},\hat{\theta}_{1},\dots,\hat{\theta}_{ k}\}
\]
is an orthonormal frame field on $(\mathbf{R}^{m+2k},g^{s})$. (Note that
$E_{\alpha}^{s}$ depends trivially on $s$ when $\alpha>m$.) Set
\[
I_{1}=\{1,\dots,m\},
\]
\[
I_{2}=\{m+1,\dots,m+k\},
\]
and
\[
I_{3}=\{m+k+1\dots,m+2k\}.
\]
\end{notation}

\begin{lemma}
\label{lemma.5.6} We use the notation of Notation \ref{rem.5.2} and Definition
\ref{def.5.3}. Let $c_{\,\,\alpha\beta}^{(s)\gamma}$ denote the structure
constants given by $[E_{\alpha}^{s},E_{\beta}^{s}]=\sum_{\gamma}
\,c_{\,\,\alpha\beta}^{(s)\gamma}E_{\gamma}^{s}$. Then:\newline {(i)} If
$\gamma\in I_{3}$ and $\alpha,\beta\in I_{1}$, then $\{c_{\,\,\alpha\beta
}^{(s)\gamma}\}$ is a homogeneous $(r,s)$-deformation of degree $-1$
.\newline {(ii)} If $\gamma\in I_{3}$, one of $\alpha,\beta$ is in $I_{1}$ and
the other is in $I_{2}$, then $\{c_{\,\,\alpha\beta}^{(s)\gamma}\}$ is a
homogeneous $(r,s)$-deformation of degree $0$.\newline {(iii)} In all other
cases $c_{\,\,\alpha\beta}^{(s)\gamma}=0$.
\end{lemma}

\begin{proof}
By Notation \ref{rem.5.2}, we have for $1\leq i,j\leq m$ and $1\leq p\leq k$,
\begin{align*}
\lbrack\hat{x}_{i},\hat{x}_{j}]  &  =\sum_{q=1}^{k}\,\left(  \frac{\partial
}{\partial x_{i}}a_{jq}^{s}-\frac{\partial}{\partial x_{j}}a_{iq}^{s}\right)
\frac{\partial}{\partial\theta_{q}}\\
&  =\sum_{q=1}^{k}\,\left(  \frac{\partial}{\partial x_{i}}a_{jq}^{s}
-\frac{\partial}{\partial x_{j}}a_{iq}^{s}\right)  r_{q}\hat{\theta}_{q}\\
\lbrack\hat{r}_{p},\hat{x}_{i}]  &  =\sum_{q=1}^{k}\, \frac{\partial}{\partial
r_{p}}a_{iq}^{s}\frac{\partial}{\partial\theta_{q}}\\
&  =\sum_{q=1}^{k}\, \frac{\partial}{\partial r_{p}}a_{iq}^{s}r_{q}
\hat{\theta}_{q} .
\end{align*}
All other brackets of vectors in our orthonormal frame are zero. The lemma now
follows from Lemma \ref{lemma.5.4} and the fact that for each $i,p$,
$\{a_{ip}^{s}\}$ is a homogeneous $(r,s)$-deformation of degree $0$.
\end{proof}

\begin{lemma}
\label{lemma.5.7} Let ${\Gamma_{\,\,\,\alpha\beta}^{(s)\gamma}}$,
$\alpha,\beta,\gamma=1,\dots,m+2k$ denote the Christoffel symbols for the
metric $g^{s}$ with respect to the frame field $\{E_{1}^{s},\dots,E_{m+2k}
^{s}\}$, i.e., $\nabla_{E_{\beta}^{s}}^{s}E_{\alpha}^{s} =\sum_{\gamma
}\,{\Gamma_{\,\,\,\alpha\beta}^{(s)\gamma}}E_{\gamma}^{s}$, where $\nabla^{s}$
is the Levi-Civita connection for $g^{s}$. Then:\newline {(i)} If two of the
indices $\alpha,\beta,\gamma$ lie in $I_{1}$ and the other lies in $I_{3}$,
then $\{{\Gamma_{\,\,\,\alpha\beta}^{(s)\gamma}}\}$ is a homogeneous
$(r,s)$-deformation of degree $-1$.\newline {(ii)} If one of the indices
$\alpha,\beta,\gamma$ lies in $I_{1}$, one lies in $I_{2}$ and one lies in
$I_{3}$, then $\{{\Gamma_{\,\,\,\alpha\beta}^{(s)\gamma}}\}$ is a homogeneous
$(r,s)$-deformation of degree 0.\newline {(iii)} If the indices $\alpha
,\beta,\gamma$ do not satisfy the conditions of either (i) or (ii), then
${\Gamma_{\,\,\,\alpha\beta}^{(s)\gamma}}=0$.\newline Thus in every case,
$\{{\Gamma_{\,\,\,\alpha\beta}^{(s)\gamma}}\}$ is a homogeneous $(r,s)$
-deformation. Moreover:\newline {(iv)} If $\{{\Gamma_{\,\,\,\alpha\beta
}^{(s)\gamma}}\}$ is a homogeneous $(r,s)$-deformation of degree $d$, then
$\{{\Gamma_{\,\,\,\alpha\beta,\delta}^{(s)\gamma}}\}$ given by ${\Gamma
_{\,\,\,\alpha\beta,\delta}^{(s)\gamma}}: =E_{\delta}({\Gamma_{\,\,\,\alpha
\beta}^{(s)\gamma}})$ is a (possibly zero) homogeneous $(r,s)$-deformation. It
is of degree $d+1$ if $\delta\in I_{2}$ and of degree $d$ otherwise.
\end{lemma}

\begin{proof}
We have
\begin{align*}
{\Gamma_{\,\,\,\alpha\beta}^{(s)\gamma}}  &  =\frac{1}{2}\{g^{s}([E_{\beta
},E_{\alpha}],E_{\gamma})+g^{s}([E_{\gamma},E_{\alpha}],E_{\beta}
)+g^{s}([E_{\gamma},E_{\beta}],E_{\alpha})\}\\
&  =\frac{1}{2}(c_{\,\,\beta\alpha}^{(s)\gamma}+c_{\,\,\gamma\alpha}
^{(s)\beta}+c_{\,\,\gamma\beta}^{(s)\alpha}).
\end{align*}
Thus this lemma follows from Lemma \ref{lemma.5.6}.
\end{proof}

\begin{proposition}
\label{prop.5.8} Let ${R_{\alpha\beta\gamma\delta}^{(s)}} =g^{s}(R^{(s)}
(E_{\gamma}^{s},E_{\delta}^{s})E_{\beta}^{s},E_{\alpha}^{s})$. Then
$\{{R_{\alpha\beta\gamma\delta}^{(s)}}\}$ is an $(r,s)$-deformation of degree
at most one. Moreover, if $\{{R_{\alpha\beta\gamma\delta}^{(s)}}\}$ has degree
one, then one of the indices $\alpha,\beta,\gamma,\delta$ lies in $I_{1}$, two
(including at least one of $\gamma$, $\delta$) lie in $I_{2}$ and one lies in
$I_{3}$.
\end{proposition}

\begin{proof}
We have
\[
{R_{\alpha\beta\gamma\delta}^{(s)}} =\Gamma_{\,\,\,\mu\gamma}^{(s)\alpha}
\Gamma_{\,\,\,\beta\delta}^{(s)\mu}-\Gamma_{\,\,\,\mu\delta}^{(s)\alpha}
\Gamma_{\,\,\,\beta\gamma}^{(s)\mu}-c_{\,\,\gamma\delta}^{(s)\mu}
\Gamma_{\,\,\,\beta\mu}^{(s)\alpha}+\Gamma_{\,\,\,\beta\delta,\gamma
}^{(s)\alpha}-\Gamma_{\,\,\,\beta\gamma,\delta}^{(s)\alpha}.
\]
By Lemmas \ref{lemma.5.4} and \ref{lemma.5.7}, each of the first three terms
belong to $(r,s)$-deformations of non-positive degree. Also by Lemma
\ref{lemma.5.7}, $\{\Gamma_{\,\,\,\beta\delta,\gamma}^{(s)\alpha}\}$ is a
homogeneous $(r,s)$-deformation. Its degree is one if $\gamma\in I_{2}$ and if
$\Gamma_{\,\,\,\beta\delta}^{(s)\alpha}$ has degree zero; otherwise its degree
is non-positive. Again applying Lemma \ref{lemma.5.7}, if $\Gamma
_{\,\,\,\beta\delta}^{(s)\alpha}$ has degree zero and is non-trivial, then
each of $I_{1}$, $I_{2}$ and $I_{3}$ contains exactly one of $\alpha$, $\beta
$, and $\delta$. The proposition now follows.
\end{proof}

\emph{Proof of Theorem \ref{thm.5.1}}: It follows from Proposition
\ref{prop.5.8}, that the integrand $\{5\tau^{2} -\Vert\operatorname{Ric}%
^{(s)}\Vert^{2}+ \Vert R^{(s)}\Vert^{2}\}$ in $\{a_{2}(g^{s})\}$ is an
$(r,s)$-deformation of degree at most $2$. We now show that it in fact has
degree 2. We have
\[
\operatorname{Ric}_{\alpha\beta}^{(s)} =\sum_{\gamma}\,R_{\alpha\gamma
\beta\gamma}^{(s)}
\]
and
\[
\Vert\operatorname{Ric}^{(s)}\Vert^{2}=\sum_{\alpha,\beta}
\,|\operatorname{Ric}_{\alpha\beta}^{(s)}|^{2}.
\]
By Proposition \ref{prop.5.8}, $\{\operatorname{Ric}_{\alpha\beta}^{(s)}\}$ is
an $(r,s)$-deformation of degree at most one. Moreover, $\{\operatorname{Ric}
_{\alpha\beta}^{(s)}\}$ has degree one only when one of $\alpha$, $\beta$ lies
in $I_{1}$ and the other in $I_{3}$. In this case, the homogeneous term of
degree one in $\operatorname{Ric}_{\alpha\beta}^{(s)}$ is equal to the
homogeneous term of degree one in $R_{\alpha\,(m+1)\,\beta\,(m+1)}
^{(s)}+\dots+R_{\alpha\,(m+k)\,\beta\,(m+k)}^{(s)}$. Consequently,
$\{\Vert\operatorname{Ric} ^{(s)}\Vert^{2}\}$ is an $(r,s)$-deformation of
degree at most two, and the homogeneous term of degree two equals
\[
2\sum_{\alpha\in I_{1}}\sum_{\beta\in I_{3}}\,|h_{1}(R_{\alpha\,(m+1)\,\beta
\,(m+1)}^{(s)})+\dots+h_{1}(R_{\alpha\,(m+k)\,\beta\,(m+k)}^{(s)})|^{2},
\]
where we are using the notation $h_{1}(\cdot)$ to denote the homogeneous term
of degree one. (The coefficient 2 is due to the symmetry when we interchange
the roles of $\alpha$ and $\beta$.) From the identity $(a_{1}+\dots a_{k}%
)^{2}\leq k(a_{1}^{2}+\dots+a_{k}^{2})$, we thus see that the homogeneous term
of degree 2 in $\Vert\operatorname{Ric} ^{(s)}\Vert^{2}$ is no bigger than
$\{2kT^{s}\}$ where
\[
T^{s}:=\sum_{\alpha\in I_{1}}\sum_{\beta\in I_{3}}\,(|h_{1}(R_{\alpha
\,(m+1)\,\beta\,(m+1)}^{(s)})|^{2}+\dots+|h_{1}(R_{\alpha\,(m+k)\,\beta
\,(m+k)}^{(s)})|^{2}).
\]
Next consider $\Vert R^{(s)}\Vert^{2} =\sum_{\alpha,\beta,\gamma,\delta
}\,\left|  R_{\alpha\beta\gamma\delta}^{(s)}\right|  ^{2}$. Note that $\{\Vert
R^{(s)}\Vert^{2}\}$ is also an $(r,s)$-deformation of degree at most two. The
homogeneous term of degree two equals the homogeneous term of degree two in
the sum of those $\left|  R_{\alpha\beta\gamma\delta}^{(s)}\right|  ^{2}$ for
which two of the indices lie in $I_{2}$, one lies in $I_{1}$ and one in
$I_{3}$. Due to the symmetries of the curvature tensor, we conclude that the
homogeneous term of degree two in $\Vert R^{(s)}\Vert^{2}$ is given by $8T(s)$.

We conclude that the homogeneous term of degree two in $\{-\Vert
\operatorname{Ric}\Vert^{2}+\Vert R^{(s)}\Vert^{2}\}$ satisfies $f^{s}%
\geq(8-2k)T^{s}\geq0$. (The latter inequality uses the hypothesis that
$k\leq3$.) Moreover, the homogeneous term of degree two in $(\tau^{s})^{2}$
equals $(\tau_{1}^{s})^{2}\geq0$, where $\tau_{1}^{s}$ is the homogeneous term
of degree one in $\tau^{s}$. By Lemma \ref{lemma.5.5} , the heat invariant
$a_{2}(g^{s})$ is a finite real linear combination of powers of $s$ with the
coefficient of $s^{2-2k}$ being given by $\int_{\mathbf{R}^{m+2k}}\,5(\tau
_{1}^{1})^{2}+f^{1}$. Thus the theorem will follow if we show that
$\int_{\mathbf{R}^{m+2k}}\,f^{1}\neq0$. Since $f^{1}\geq0$, we have
$\int_{\mathbf{R}^{m+2k}}\,f^{1}=0$ only if $f^{1}\equiv0$.

Choose $\alpha,\beta,\gamma,\delta$ as follows: Choose $i\in\{1,\dots,m\}$ and
$p\in\{1,\dots,k\}$ so that the linear transformation ${\mathbf{R}}^{m}
\rightarrow\mathfrak{z}$ given by $x\rightarrow\langle\lbrack x,e_{i}
],Z_{p}\rangle$ is not identically zero. Without loss of generality, we assume
$p=1$. Let $\beta=\delta=m+2$, let $\alpha=m+k+1$ and let $\gamma=i$. Then by
Notation \ref{rem.5.2} and the proofs of Lemma \ref{lemma.5.6}, Lemma
\ref{lemma.5.7}, and Proposition \ref{prop.5.8}, the homogeneous term of
degree one in $R_{\alpha\beta\gamma\delta}^{1}$ is given by
\[
\frac{\partial^{2}}{\partial r_{2}^{2}}a_{i1}^{1}\,r_{1} =4r_{2}^{2}\phi
_{22}(\Vert x\Vert^{2},\Vert r\Vert^{2})\langle\lbrack x,e_{i}],Z_{1}\rangle
r_{1}
\]
(viewed as a function on ${\mathbf{R}}^{m+2k}$ depending trivially on $\theta
$). This function cannot be identically zero since the smooth cut-off function
$\varphi$ cannot be linear in either variable. It follows that $f^{1}$ is not
identically zero. This completes the proof of Theorem \ref{thm.5.1}. $\square$

\end{document}